\UseRawInputEncoding

\documentclass[12pt]{amsart}
\usepackage{times}
\usepackage{mathptmx}
\makeatother

\theoremstyle{definition}

\def\fnum{equation} 
\newtheorem{Thm}[\fnum]{Theorem}
\newtheorem{Cor}[\fnum]{Corollary}

\newtheorem{Lem}[\fnum]{Lemma}

\newtheorem{Def}[\fnum]{Definition}
\newtheorem{Exa}[\fnum]{Example}

\newtheorem{Pro}[\fnum]{Proposition}

\numberwithin{equation}{section}

 \newcommand{\E}{\ensuremath{\mathbb{E}}}

 \newcommand{\ba}{\begin{align*}}
 \newcommand{\ea}{\end{align*}}

\title{A NEW DEFINITION OF THE DIMENSION OF GRAPHS }

\author{X.Zhao}%
\address{}

\thanks{}

\email{}

\begin{document}

\maketitle

 \begin{abstract}
Enlighted by the recent work of Hao Huang on sensitivity conjecture \cite{2}, we propose a new definition of the dimension of graphs and establish a relationship between the chromatic number and the dimension.
  \end{abstract}
\section{introduction}

Dimension is a well-established and crucial notion in both geometry and algebra, which reflects the degrees of freedom and plays an essential role in classifying objects. It reveals, especially in geometry, the local picture of the subject. Since they are all in the continuous setting, it is natural to ask what happens when moving to the discrete case. In this paper, I would like to establish a new definition of the dimension of graphs and understand how it plays with other properties of graphs.

As mentioned by Erd$\ddot{o}$̈s, Harary and Tutte \cite{1}, there is a concept of dimension of graphs defined according to the minimum dimension of Euclidean space such that the graph can be put into keeping all edges of length 1 and vertices distinct. It reflects the complexity of graphs in the sense of rigidity, which is more a global property instead of telling what happens locally. Moreover, since the main concern is on the distance between vertices and the intersections between edges are permitted, some graphs may have extremely low dimension benefiting from the regularity of their structure. For example, $Q_n$  ($n\geq2$)  are all of dimensions 2 in that sense, which is kind of counter-intuitive as they are certainly different locally. 

Motivated by the resent work of Huang \cite{2} on sensitivity conjecture, which reveals that for any $2^n+1$ points in $Q_n$, there must be a point which has degree not less than $\sqrt{n}$. It inspires me a new definition of dimension of graphs, which could reflect the more local difference of graphs than the former one. On the other hand, I will show that in some sense this dimension controls the former one from above.
\section{New Definition of Dimension}
Now let’s come to the definition part. Given a graph $G=(V,E)$, where V stands for the vertices and E stands for the edges. Let $G'=(V',E')$ be a sub-graph of G, then we first define the sub-dimension of G' as 
\begin{Def}
\begin{align*}
subdim(G')=\min_{\mbox{\tiny$\begin{array}{c}
G''=(V'',E'')\subset G'\\|V''|\geq |\frac{|G'|}{2}+1\\E''\,\,is\,\,maximal
\end{array}$}} \max_{x\in V'' q}deg(x)
\end{align*}
\end{Def}           
Or more concisely with the notation $\delta(G'')=\max_{x\in V''}⁡deg⁡(x)$
\begin{align*}
subdim(G')=\min_{\mbox{\tiny$\begin{array}{c}
G''=(V'',E'')\subset G'\\|V''|\geq |\frac{|G'|}{2}+1\\E''\,\,is\,\,maximal
\end{array}$}}\delta(G'')
\end{align*}
And finally, we define
\begin{Def}
\begin{align*}
dim(G')=\max_{\mbox{\tiny$\begin{array}{c}
G'=(V',E')\subset G'\\E''\,\,is\,\,maximal
\end{array}$}} subdim(G')
\end{align*}
\end{Def} 
Here E'' is maximal means that it contains all the edges inside  E' whose vertices belong to V''.

The reason why we define in this way is the following. Since we want the definition reflects the difference between $Q_n$, it is very natural to introduce the dimension of sub-dimension based on the work of Huang, which show that they are $\sqrt{n}$ for $Q_n$.

The primary reason of taking the definition as the maximum among all sub-dimension of the sub-graphs is that we want the concept of dimension to reflect the local behavior of the graph, which indicates that it should be at least non-decreasing with respect to the inclusion relation. Another reason is that sub-dimension is sometimes degenerate which is certainly not wanted, while the dimension is always non-zero for graphs with edges from the monotonicity. 

It is deserved to be pointed out that when calculating the sub-dimension, we take sub-graphs which have more than $[\frac{|G' |}{2}]+1$ points. Here we take the floor function.

To make it more clearly, here I would like to give some simple examples.

\begin{Exa}
For a linear graph with n ($n
\geq 4$) points, the dimension is 1 as when we take more than [n/2]+1 points from it, there must be two of them adjacent from the drawer principle. On the other hand we can take 1, 2, 4, 5, … $3i+1$, $3i+2$… which clearly realizes the lower bound 1. Since the sub-dimension are all 1, the dimension is 1 from the definition. 
\end{Exa}
  
\begin{Exa}
     For loops with n ($n\geq 5$) points, similar as before, we can use the drawer principle and the specific construction in example to conclude that they are all of dimension one. It is reasonable because locally a loop looks exactly as linear graph while the difference is shown only when we move to the global picture, which is not supposed to be reflected by its dimension. \end{Exa}
  
\begin{Exa}

     For complete graph $K_n$, no matter taking any $[\frac{n}{2}]+1$ points from it with maximal edge set, it would form a $K_([\frac{n}{2}]+1)$, which has degree $[\frac{n}{2}]$. So the sub-dimension of $K_n$ is $[\frac{n}{2}]$. Since any sub-graph with maximal edges of $K_n$ is a complete graph, we know that its sub-dimension is at most $[\frac{n}{2}]$. Thus the dimension of $K_n$ is $[\frac{n}{2}]$.
\end{Exa}
  
\begin{Exa}
For complete bipartite graph $K_{(m,n)}$  ($m\leq n$), if $n\neq m$, then its sub-dimension is 0 since we can take the points in the side corresponding to n. If $n=m$, then the best strategy is to take $[\frac{m}{2}]+1$ and $m-[\frac{m}{2}]$ respectively and the sub-dimension is $[\frac{m}{2}]+1$. Since any sub-graph of $K_{(m,n)}$ with maximal edge set is of the form $K_{(a,b)}$  $(a\leq m,b\leq n,a\leq b)$. So the dimension of $K_{(m,n)}$  ($m\leq n$) is $[\frac{m}{2}]+1$ 
\end{Exa}\section{property of the dimension}
After considering these examples, it might be natural to come back to our motivation, that is, the $Q_n$. Huang’s work \cite{2} is equal to say that the sub-dimension of $Q_n$ is $[\sqrt{n}]$. As we have defined the dimension as the maximum of sub-dimension of sub-graphs with maximal edge set, it is not clear at this point whether $Q_n$’s dimension is $[\sqrt{n}]$ as a sub-graph of it might be weird. It is natural, however, to guess that the dimension of $Q_n$ is also  $[\sqrt{n}]$, while it is indeed the case using the symmetric structure of $Q_n$.

\begin{Thm}
\begin{align*}dim(Q_n)=[\sqrt{n}]\end{align*}\label{thm1}
\end{Thm}              

\begin{proof}
From Huang’s work \cite{2}, we know that the sub-dimension of $Q_n$ is at least $[\sqrt{n}]$. To verify that the lower bounded is achieved, we would need to introduce a construction given by Chung \cite{3}.
\end{proof}
\begin{Lem}
    There is a sub-graph $G_n$ of $Q_n$ with $2^(n-1)+1$ vertices and maximal edges such that the degree of $G_n$ is $[\sqrt{n}]$\label{lem1}
\end{Lem}
    
Using this lemma, we are now ready to prove that the dimension of $Q_n$   is $[\sqrt{n}]$. It is equal to say that any sub-graph of $Q_n$ with maximal edge set has sub-dimension less or equal to $[\sqrt{n}]$. Take a sub-graph $G'=(V',E' )\subset Q_n$ with maximal edge set $E'$, which means that it contains all the edges of $Q_n$ between vertices in $V'$. It now suffices to prove that the sub-dimension of $G'$ is equal or less than $[\sqrt{n}]$. Here we use the symmetry of $Q_n$. Denote the points in $Q_n$ as n tuples $(a_1,a_2,…,a_(n-1),a_n)$ where $a_i\in\{0,1\}$. Then the $G_n$ mentioned in the lemma above is a subset of $2^(n-1)+1$ vertices while the edge set is maximal. Then for each $(a_1,a_2,…,a_(n-1),a_n)\in Q_n$, we define
\begin{align*}
    G_{(a_1,a_2,…,a_(n-1),a_n )}=\{x+(a_1,a_2,…,a_(n-1),a_n ):  x\in G_n\}
\end{align*}
where addition is defined in $Z_2$. So we have $2^n$ such sets while each has $2^(n-1)+1$ points.

Consider the intersection of them with $G'$, the sub-graph we fixed before. Here we will use the probability type argument. Since each $ G_{(a_1,a_2,…,a_(n-1),a_n )}$ contains $2^(n-1)+1$ points and the symmetry of $Q_n$, each point in $G'$ shows up exactly $2^(n-1)+1$ times in ${ G_{(a_1,a_2,…,a_(n-1),a_n )}: (a_1,a_2,…,a_(n-1),a_n )\in Q_n}$. So we know that 
\begin{align*}
    \sum_{(a_1,a_2,…,a_(n-1),a_n )\in Q_n}\#G_{(a_1,a_2,…,a_(n-1),a_n )}\cap G'=(2^{n-1}+1)|G'|
\end{align*}

Since the left-hand side is a sum of $2^n$ numbers, there must be at least one which is no less than $\frac{2^(n-1)+1}{2^n} |G' |$. In other words, there is a $(a_1,a_2,…,a_(n-1),a_n )\in Q_n$ such that \begin{align*}
\#G_{(a_1,a_2,…,a_(n-1),a_n )}\cap G'\geq \frac{1}{2}|G'|+1
\end{align*}

From the property of $G_{(a_1,a_2,…,a_(n-1),a_n )}$, which is the same as $G_n$ in the sense of translation, we know that the degree of $G_{(a_1,a_2,…,a_(n-1),a_n )}$ is $[\sqrt{n}]$ as shown by the lemma.

Since $\#G_{(a_1,a_2,…,a_(n-1),a_n )}\cap G'\geq \frac{1}{2}|G'|+1$ and $G_{(a_1,a_2,…,a_(n-1),a_n )}\cap G'$ is a sub-graph of both $G'$ and $G_{(a_1,a_2,…,a_(n-1),a_n )}$ with maximal edge set, the degree of $G_{(a_1,a_2,…,a_(n-1),a_n )}\cap G'$ is controlled by that of $G_{(a_1,a_2,…,a_(n-1),a_n )}$, which is $[\sqrt{n}]$. So the sub-dimension of $G'$ is not more than $[\sqrt{n}$. Combined with the construction of $G'$, we concluded that the dimension of $Q_n$ is $[\sqrt{n}]$.       

With exactly the same argument, we can show that 
\begin{Pro}
For any Cayley graph, its dimension equals to its sub-dimension\label{pro1}
\end{Pro}

\begin{proof}
Use the fact that G has group structure to construct the $G_a$ as that in Theorem \ref{thm1} 
\end{proof}
This theorem shows that the definition of dimension is reasonable and realizes our motivation, which indeed reveals the local difference among $Q_n$.

Now since we are in the world of graphs, it would be natural to ask what the interaction between the dimension and other concepts in graph theory is. Since the concept is completely new, I believe there would be many results waiting to be found. Here we will give some fundamental results which could be enlightening.

Firstly we would like to show that the dimension controls the chromatic number $\chi$, which is the minimum number of colors needed to color the vertices such that no two adjacent vertices having the same color. 
\begin{Thm}
\begin{align*}
    \chi(G)\leq (\dim(G)+1)\cdot |log_2(|G|)|
\end{align*}\label{thm2}
\end{Thm}              

\begin{proof}
To prove this theorem, we need to characterize the chromatic number in the language of degrees, which is stated in the following lemma

\begin{Lem}
    If G is of chromatic number n, then there is a sub-graph $G'$  of G such that each point in $G'$ has degree not less than $n-1$
\end{Lem}
    
\begin{proof}
We consider the critical sub-graph, which means that taking out any vertex will make the chromatic number decrease 1. It is deserved to be pointed out that taking any point out of a graph will either keep the chromatic number or let it reduce 1. Since we can keep taking out vertices that are not relevant to the chromatic number, we will ultimately reach a critical sub-graph $G'$ of G with chromatic n.   

Since it is critical, we know that taking out any point will make its chromatic number to be $n-1$. Suppose there is indeed a point x in $G'$ with degree less than $n-1$, then we can first take x out and color the rest points. After that we take back x and choose a suitable color among those $n-1$ ones, which is possible since x has degree less than $n-1$ and forbids at most $n-2$ colors. This construction contradicts with the fact that $G'$ has chromatic number n. So we know the assumption that x has degree less than $n-1$ is wrong, which proves the lemma.\end{proof} 

Now we come back to the proof of the Theorem \ref{thm2}. From the definition of dimension, we know that there is a sub-graph $G'$ of G with $[\frac{|G|}{2}]+1$ point and has degree not more than dim(G), which means that each point of $G'$ has degree less or equal to dim(G). From the lemma we know that the chromatic number of $G'$ is at most $dim(G)+1$. After that we delete all the points of $G'$ from G and continue what we did before. Since the dimension is monotonic, dim(G) is not less than dim($G'$). Continue the process until we get a point at last. 

In this process we get split the vertex set of G into at most $[log_2 |G|]$ distinct groups, each having chromatic number at most $dim(G)+1$. So we can now color them in the way that vertices in the different groups has completely different colors, while inside the same group we know their chromatic number is not more than $dim(G)+1$ from the argument above. Thus the chromatic number of G is at most $(dim⁡(G)+1)\cdot[log_2(|G|)]$.

In the proof above we apply the strategy that split the whole graph into several distinct sub-graphs and treat them separately. One might try to find a finer estimate with more consideration on the interaction between the sub-graphs mentioned above, which could be tricky and needs more information on the structure of the graph.

Recall that there is a former definition of graph proposed by Erd$\ddot{o} ̈$s, Harary and Tutte i n\cite{1}, which we may denote it as $\dim_{EHT}$

Since it is not good at revealing the local structure of the graph, we proposed the new definition mention above. Now it would be natural to think about the relationship between those two concepts of dimension. Now we can show that in some sense  $\dim_{EHT}$ is controlled by dim from above.

\begin{Cor}\begin{align*}
    \dim_{EHT}(G)\leq 2(dim⁡(G)+1)\cdot[log_2(|G|)
\end{align*}\label{cor1}
 
\end{Cor}

\begin{proof}
Here we introduce a Lemma proposed in \cite{1}

\begin{Lem}(\cite{1})
 \begin{align*}
     \dim_{EHT}(G)\leq 2\chi(G)
 \end{align*}\label{lem3}
\end{Lem}\begin{proof}
Firstly, we reserve two dimensions for each color, which means that we need $2\chi(G)$ dimensions in total. Set the coordinates such that points are all in the $\frac{1}{\sqrt{2}}$-circle and have non-zero coordinates only when restricting to the two dimensions corresponding to their color. For example, those with the color 1 has coordinates like $(a,b,0,0,…,0)$, where $a^2+b^2=\frac{1}{2}$. Since any two adjacent points must have different colors, we know that the distance between them must be 1 for Lemma \ref{lem3}\end{proof}   
Now the corollary 1 is immediately gotten from Theorem \ref{thm2} and  Lemma \ref{lem3}, because $  \dim_{EHT}(G)\leq 2\chi(G)\leq 2(dim⁡(G)+1)\cdot[log_2(|G|)]$\end{proof}   \end{proof}

\end{document}